\documentclass[letterpaper, 10pt, conference]{ieeeconf}
\IEEEoverridecommandlockouts
\makeatletter
	\let\NAT@parse\undefined
\makeatother

\PassOptionsToPackage{noend}{algpseudocode}
\PassOptionsToPackage{nameinlink,capitalize}{cleveref}
\PassOptionsToPackage{dvipsnames}{xcolor}
\let\labelindent\relax
\usepackage[nocompress]{cite}
\usepackage[
	noload={todonotes,tikz,pgfplots},
]{myPreamble}

\renewcommand{\ifcaseslineend}{}
\newcommand{\email}[1]{\href{mailto:#1}{#1}}

\crefname{figure}{Figure}{Figures}
\crefname{ALG@line}{step}{steps}
\newcommand{\secref}[1]{\hyperref[#1]{\S\ref*{#1}}}
\newcommand{\figref}[1]{\hyperref[#1]{Fig. \ref*{#1}}}

\usepackage{pifont}

\newcommand{\figname}{}
\newenvironment{addHFigTile}[2][]{%
	\begin{subfigure}[t]{.5\linewidth}%
	\captionsetup{width=.95\linewidth}%
	\vspace*{0pt}%
	\fbox{\includegraphics[width=.97\linewidth,#1]{Pics/\figname#2}}%
}{%
	\end{subfigure}%
	\ignorespacesafterend
}

\newenvironment{addVFigTile}[2][]{%
	\begin{subfigure}[t]{.25\linewidth}%
		\captionsetup{width=.95\linewidth}%
		\vspace*{0pt}%
		\fbox{\includegraphics[angle=90,origin=c,width=.94\linewidth,#1]{Pics/\figname#2}}%
}{%
	\end{subfigure}%
	\ignorespacesafterend
}

\title{\LARGE\bf
	Flock navigation with dynamic hierarchy and subjective weights
	\texorpdfstring{\\}{}using nonlinear MPC%
}
\author{\Large
	Aneek Nag\qquad
	Shuo Huang\qquad
	Andreas Themelis\qquad
	Kaoru Yamamoto\thanks{%
		Kyushu University, Department of Electrical Engineering,
		744 Motooka, Nishi-ku, 819-0395 Fukuoka, Japan.
		This work is supported by the Japan Society for the Promotion of Science (JSPS) KAKENHI grants n. JP19H02161 JP20K14766 and JP21K17710.
		The first author is supported by the Japanese government (MEXT) scholarship.%
		\sf
		\newline
		\email{nag.aneek.364@s.kyushu-u.ac.jp},
		\email{huangshuo0415@outlook.com},%
		\newline
		\email{andreas.themelis@ees.kyushu-u.ac.jp},
		\email{yamamoto@ees.kyushu-u.ac.jp}%
	}%
}

\begin{document}

	\maketitle

	\begin{abstract}
		We propose a model predictive control (MPC) based approach to a flock control problem with obstacle avoidance capability in a leader-follower framework, utilizing the future trajectory prediction computed by each agent.
We employ the traditional Reynolds' flocking rules (cohesion, separation, and alignment) as a basis, and tailor the model to fit a navigation (as opposed to formation) purpose.
In particular, we introduce several concepts such as the credibility and the importance of the gathered information from neighbors, and dynamic trade-offs between references.
They are based on the observations that near-future predictions are more reliable, agents closer to leaders are implicit carriers of more educated information, and the predominance of either cohesion or alignment is dictated by the distance between the agent and its neighbors.
These features are incorporated in the MPC formulation, and their advantages are discussed through numerical simulations.

	\end{abstract}

	\section{Introduction}
		Nature often exhibits collective behavior, such as in a school of fish or a flock of birds, in which each member seemingly follows simple principles using only local information and still achieves a global task.
This is an attractive feature when dealing with the control of a large number of autonomous vehicles, especially if each one has limited computational capability.
In fact, in such a system each agent acts based solely on the information it receives (or detects) from other agents in its near proximity, without the need to be aware of the state of the entire flock.

In 1987, Reynolds introduced a model where each agent, referred to as a \emph{boid}, obeys the following three simple rules to achieve \emph{flocking} \cite{reynolds1987flocks}:
\emph{cohesion}, promoting proximity among boids by biasing each one towards the local center of mass,
\emph{separation}, which imposes a minimum safety distance so as to prevent collisions,
and \emph{alignment}, encouraging agreement of velocity among neighboring boids.
Because of its simplicity and scalability, Reynolds' flocking has received much attention ever since, and many attempts have been made to employ this model for the joint control of large numbers of vehicles and mobile robots \cite{beaver2021overview}. 
Each boid's behavior can be captured by imposing an artificial potential field representing attractive/repulsive forces and velocity consensus to the neighbors; see, e.g., \cite{tanner2003stableI,tanner2003stableII,olfatisaber2006flocking}.
However, potential-field based methods tend to exhibit oscillatory behavior as each agent reacts pretty much passively to the current state of the surroundings only.
Instead, we may consider a ``smarter'' flock, in which each boid accounts for future stages and acts accordingly.
This is precisely the idea behind MPC, a closed-loop control technique that provides optimal inputs by computing solutions of a finite-horizon problem at each sampling time, and then implementing the control law in accordance to a receding horizon strategy.
In doing so, it can account for complex problem formulations, explicitly handles state constraints, and in considering long prediction horizons it favors smoother control actions.
The applicability of MPC is however tightly linked to the efficiency of the optimization method employed to solve each problem, as well as the capability of the computing platform at hand.
Nevertheless, the outstanding computing advancements and the development of efficient tailored algorithms \cite{ohtsuka2004continuation,patrinos2014accelerated,stella2017simple,themelis2018proximal,sopasakis2020open} have made a reliable employment of this control technology a reality, even on embedded devices.

Attracted by the decentralized nature of flocks and the potential of MPC, several works for flock control using MPC strategies have been proposed in the last decade in which agents cooperatively exchange information to minimize the deviation from Reynolds' rules throughout long prediction horizons \cite{zhan2011flocking,zhan2013flocking,zhang2015model,zhang2016model,zhou2017distributed,lyu2021multivehicle,soria2022distributed}.

\subsection{Goals and proposed methodology}%
	Inspired by this trend, in this paper we aim at investigating how (suitable modifications of) flocking rules can conveniently be exploited for \emph{navigating} unknown environments towards a desired destination.
	Specifically, we consider the scenario in which the destination is unknown to some (possibly all but one) agents, yet the enforcement of a boids-type behavior allows the whole swarm information to be passively passed upon each single agent.
	In this regard, we refine the purpose of merely achieving \emph{flocking}, as in forming and maintaining a stable configuration, and focus instead on a distributed \emph{navigation} control.
	To seamlessly achieve this task, we introduce several factors to the proposed MPC formulation and to the flock modeling.

	\subsubsection{Processing the neighbors' predicted information}%
		Similarly to and inspired by \cite{yang2008formation,lindqvist2021scalable,soria2022distributed}, we propose an MPC-based approach where at each sampling time neighboring agents exchange their predicted configurations to facilitate the achievement of a cooperative task.
		To utilize such information wisely, we introduce a discount factor to reduce the influence of distant-future prediction as it is less accurate than the near-future one.
		Possible MPC horizons mismatch and delayed information are handled suitably.
	
	\subsubsection{Flock modeling}%
		As a means to expedite information flow, we design a real-time adjustment of the flocking rules that, while retaining a fully decentralized nature, biases the instantaneous cohesion and alignment configurations to better reflect the navigation objective of the entire swarm.
		In particular, we introduce the following three adjustments to the original boids model: \emph{(i)} time-varying leader-follower graph-distance-based hierarchy, \emph{(ii)} subjective neighbor weighting and flocking references, and \emph{(iii)} cohesion/alignment dynamic trade-off.
		We note that hierarchy structures in flocking have been studied in \cite{shen2008cucker} for strict hierarchy (a lower-ranked agent does not affect a higher-ranked agent at all) and in \cite{shen2008cucker,jia2019modelling} for time-invariant hierarchy.
		We instead consider a non-strict and time-variant hierarchy that each agent can compute in a decentralized manner at negligible cost.

	The paper focuses on the flock modeling, and assumes that the problem instances can be solved to \emph{local} optimality in real time.
	While not tied to a particular optimization method, we point out the OpEn solver \cite{sopasakis2020open} as a canditate well suited to our purposes, which we use 
	in the simulations.

\subsection{Paper organization}
	\Cref{sec:setting} describes in detail the investigated setting and defines the MPC problems that each boid agent needs to solve.
	Technical considerations on the information exchange are also discussed.
	\Cref{sec:flocking} presents the proposed dynamic modifications of the flocking rules, elaborating on their \emph{hierarchical} and \emph{subjective} nature.
	In \cref{sec:simulations} we showcase the advantages of the proposed modifications with numerical simulations.
	\Cref{sec:conclusions} concludes the paper.

	\section{Problem setting}\label{sec:setting}%
		Our goal is the design of a distributed control strategy for a fleet of \(N\) agents, divided into \emph{leaders}, that know (or are in charge of planning) the trajectory to reach a desired destination, and \emph{followers}, which are instead only aware of their immediate surroundings.
Specifically, we aim at controlling these latter ones under the assumption that the leaders dispatch their near future trajectory and velocity predictions to the followers, yet without necessarily taking into account infomation received from them.
Desirable features of the control action include minimizing the time to reach destination, as well as the number of required leaders and the amount of constraints on their motion for ensuring that no follower is left behind.
The revisited flocking rules developed in \cref{sec:flocking} are designed in this perspective.

\subsection{Notational conventions}
	We bold-face variables indexed along the time steps spanning the prediction horizon, e.g.,
	\(\bm u_i^t=u_i^{t|t}\ldots u_i^{t+T_i-1|t}\) shall denote a sequence of inputs of agent \(i\) and similarly
	\(\bm x_i^t\) will represent a sequence of states througout the prediction horizon \(T_i\) of agent \(i\).
	At every time instant \(t\), agent \(i\) is aware of a set of obstacles \(\mathcal O_i^t\) in its proximity and receives a sequence of predicted ouputs \(\bm y_{j|i}^t\) from each of its current in-neighbors \(j\in\mathcal N_i^t\setminus\set{i}\).
	Outputs comprise position \(\bm p_{j|i}^t\) and velocity \(\bm v_{j|i}^t\); the indexing \(j|i\) accounts for possible mismatches in the same information collected by different agents.

\subsection{MPC formulation}
	The discrete time step is considered normalized as $\Delta t=1$ time units for notational simplicity.
	At each time step \(t\), agent \(i\) needs to solve the following problem
	\begin{subequations}\label{eq:P}
		\begin{align}
			\minimize_{
				\bm u_i^t\in\mathcal U_i\times\dots\times\mathcal U_i
			}{}~
		&
			J(\bm u_i^t)
		\\
			\stt{}~
		&
			\begin{cases}[l@{~}l@{}]
				\left.
					\begin{array}{@{}rcl}
						x_i^{t+k+1|t} &{}={}& f_i(x_i^{t+k|t},u_i^{t+k|t})\\
						x_i^{t+k+1|t} &{}\in{}& \mathcal X_i\\
						\hphantom{d_{j|i}^{t+k+1|t}}
						\mathllap{p_i^{t+k+1|t}} &{}\notin{}& O,\ O\in\mathcal O_i^t
					\end{array}
				~\quad\right]
				&
				k<T_i
			\\
				\left.d_{j|i}^{t+k+1|t}\geq d_{i,\rm sep},\ j\in\mathcal N_i^{t+k+1|t}\right.
				&
				k<T_{i,\rm sep}
			\end{cases}
		\end{align}
	\end{subequations}
	where \(\func{f_i}{\R^{n_{x,i}}\times\R^{n_{u,i}}}{\R^{n_{x,i}}}\) are agent \(i\)'s discrete state dynamics, \(\mathcal U_i\subseteq\R^{n_{u,i}}\) and \(\mathcal X_i\subseteq\R^{n_{x,i}}\) are input and state constraint sets, and \(T_{i,\rm sep}<T_i\) and \(d_{i,\rm sep}>0\) are modeling parameters that account for the separation rule (see \cref{sec:discount}).
	Furthermore, at every time step \(\tau\in[t+1,t+T_i]\),
	\[
		d_{j|i}^{\tau|t}
	{}\coloneqq{}
		\bigl\|p_i^{\tau|t}-p_{j|i}^{\tau-1|t}\bigr\|^2
	\]
	is the predicted (square) distance from neighbor \(j\), \(\bar y_i^{\tau|t}\) is a reference output comprising of vectors
	\begin{equation}\label{eq:reference}
		\bar p_i^{\tau|t}
	{}\coloneqq{}
		\sum_{j\in\mathcal N_i^{\tau|t}}w_{j|i}^{\tau|t}\cdot p_{j|i}^{\tau|t}
	\quad\text{and}\quad
		\bar v_i^{\tau|t}
	{}\coloneqq{}
		\sum_{j\in\mathcal N_i^{\tau|t}}\omega_{j|i}^{\tau|t}\cdot v_{j|i}^{\tau|t},
	\end{equation}
	namely weighted averages of positions and velocities among the available neighbors \(\mathcal N_i^{\tau|t}\) with (positive, summing to one) cohesion and alignment weights \(w_{j|i}^{\tau|t}\) and \(\omega_{j|i}^{\tau|t}\) (cf. \cref{sec:weights}), and
	\[
		J(\bm u_i^t)
	{}={}
		\bigl\|
			\bm u_i^t
		\bigr\|_{R_i}^2
		{}+{}\!
		\sum*_{k=0}^{T_i-1}\!\gamma_i^k
		\bigl\|
			y_i^{t+k+1|t}
			{}-{}
			\bar y_i^{t+k|t}
		\bigr\|_{Q_i^t}^2
		{}+{}
		\rho_{i,\rm sep}\!
		\sum*_{k>T_{i,\rm sep}}\!\gamma_i^k{
			\sum*_{j\in\mathcal N_i^{t+k|t}\mathrlap{{}\setminus\set i}}{
				d_{j|i}^{t+k|t}
			}
		}
	\]
	is the objective.
	Here, \(\gamma_i\in(0,1]\) is a discount factor, \(Q_i^t\in\R^{4\times 4}\) and \(R_i\in\R^{n_{u,i}T_i\times n_{u,i}T_i}\) are (positive semidefinite, diagonal) weight matrices, and \(\rho_{i,\rm sep}\geq0\) is a penalty parameter for soft-con\-strain\-ing (or neglecting, if chosen null) separation at late stages.
	The overall information exchange and data processing is synopsized in \cref{alg:workflow}.

	\begin{algorithm}[htb]
		\caption{Conceptual MPC workflow}
		\label{alg:workflow}
		\begin{algorithmic}[1]
\item[At the beginning of stage \(t\), boid agent \(i\)]
\State\label{state:detect}%
	{\sc detects} obstacles \(\mathcal O_i^t\) and in-neighbors \(\mathcal N_i^t\) (including itself) collecting their last
	\begin{itemize}
	\item
		hierarchy level \(\pi_j^{t-1}\) and
	\item
		time-stamped predicted outputs \(\bm y_{j|i}^{t-1}\), \(j\in\mathcal N_i^t\)
	\end{itemize}
\State
	{\sc shifts} neighbors outputs based on discrepancy with own time-stamp and/or MPC horizon, and refines future neighbors \(\mathcal N_i^{t+k|t}\), \(k=1,\dots,T_i-1\), as in \secref{sec:delays}
\State
	estimates own {\sc hierarchy level} \(\pi_i^t\) as in \eqref{eq:hierarchy}
\State
	assigns {\sc individual weights} \(w_{j|i}^{t+k|t}\) and \(\omega_{j|i}^{t+k|t}\) as in \secref{sec:weights}, \(k=0,\dots,T_i-1\), \(j\in\mathcal N_i^{t+k|t}\), thus creating
	\begin{itemize}
	\item
		reference positions \(\bar p_i^{t+k|t}\) and
	\item
		reference velocities \(\bar v_i^{t+k|t}\) as in \eqref{eq:reference}
	\end{itemize}
\State
	selects {\sc cohesion/alignment trade-off weight} \(Q_i^t\) as in \eqref{eq:Q}
\State\label{state:solve}%
	{\sc computes optimal inputs} \(\bm u_i^t\) by solving problem \eqref{eq:P}, and store predicted outputs \(\bm y_i^t\)
\State
	{\sc actuates} \(u_i^{t|t}\), \(t\gets t+1\), and {\sc restarts} from \cref{state:detect}
\end{algorithmic}

	\end{algorithm}
	
\subsection{Horizons mismatch and time delays}\label{sec:delays}%
	The dependency on \(\tau\) for the set of neighboring agents \(\mathcal N_i^{\tau|t}\) accounts for the possibility of different prediction horizons \(T_i\) among neighbors and/or delayed information.
	The rationale is as follows.
	If agent \(i\) receives information from an agent \(j\in\mathcal N_i^t\) which has a longer prediction horizon \(T_j>T_i\), then any prediction from step \(T_i\) onward can be discarded.
	In case \(T_j<T_i\), instead, agent \(j\) can be removed from the \emph{virtual} set of neighbors \(\mathcal N_i^{\tau|t}\) for \(\tau\geq T_j\), since no information is anymore available after that time.

	The strategy of \emph{virtually} removing neighbors in the prediction can also be adopted to account for delayed information.
	If agent \(i\) receives information from agent \(j\) with a delay of \(\delta\) time steps, then the received output sequence is \(\bm y_{j|i}^{t-\delta}\), of which only the last \(T_i-\delta\) elements will be taken into account, leading to \(\bm y_{j|i}^{t}=y_{j|i}^{t|t-\delta},\dots,y_{j|i}^{t-\delta+T_i|t-\delta}\), at which point any shortcomings can be accounted for by suitably removing agent \(j\) from late virtual neighborhoods.
	
	Note that, by adopting the convention that \(i\in\mathcal N_i^{\tau|t}\) for all \(\tau\), in absence of external information the average \(\bar y_i^{\tau|t}\) reduces to agent \(i\)'s previously predicted output, so that the term
	\(
		\|
			y_i^{t+k|t}
			{}-{}
			\bar y_i^{t+k|t}
		\|_{Q_i^t}^2
	\)
	in the cost function promotes consistency with the information that was dispatched to the neighbors.

\subsection{Discount factor and late-stage soft constraints}\label{sec:discount}%
	In accordance with the receding horizon principle of MPC, although a sequence of \(T_i\) optimal inputs is computed every time it is only the first one that is actuated.
	This strategy follows the logic that predictions on late stages are much less accurate and more susceptible to environmental changes, which is why the problem is continuously updated to take new information into account.
	Following this spirit and as a way to prevent too much conservatism, the separation rule (i.e., collision avoidance) is enforced only on early time instants, and is otherwise only discouraged as a soft penalty in the cost function for later stages with a quadratic term \(\rho_{i,\rm sep}d_{j|i}^{\tau|t}\).
	For the same reason, a discount factor \(\gamma_i^k\), with \(\gamma_i\in(0,1]\) weighs the \(k\)-th per-stage cost in the objective so as to prioritize early stages in the decision process and give less importance to longer term estimates which are likely less accurate.
	Avoidance of detected obstacles, instead, is hard-constrained for the entire prediction horizon, owing to the fact that obstacles are assumed static and thus not subject to (unknown) changes in position.

	\section{Flocking rules revisited}\label{sec:flocking}%
		In this section we propose some adaptions of the classical flocking rules that are better suited for navigating in an unknown environment.
The advantages of each of them will be showcased with simulations in the following section.

\subsection{Leader-follower graph-distance hierarchy}\label{sec:leader}%
	Because of their predominant role and more accurate and longer term spatial information, it is reasonable to give more importance to data received from leaders.
	For the same reason, followers that are directly connected to a leader implicitly benefit for more educated estimates than those which are not, and their dispatched data should likewise be given more weight.
	As a way to enforce this chain of priority, we may endow the connectivity graph of a weight on the (directed) edges based on the graph distance from the closest leader.
	Although each follower is only aware of a limited neighborhood and this quantity may thus not be readily available, it is nevertheless easy to well estimate it through the information exchanges.
	Each leader \(\ell\) is assigned a hierarchy level \(\pi_\ell^t\equiv0\), constant over time, while each follower \(i\), starting from a predefined upper bound \(\pi_i^0=\overline\pi\in\N_{\geq1}\), updates its value over time as
	\begin{equation}\label{eq:hierarchy}
		\pi_i^t=\min\set{\overline\pi,\,1+\min_{i\neq j\in\mathcal N_i^t}\pi_j^{t-1}}.
	\end{equation}
	The cap \(\overline\pi\) is a model parameter, possibly different among different agents and whose role is merely to prevent arbitrarily low weights.
	Although the rule is only based on local information, it is easy to see that it accurately estimates the graph distance from the set of leader agent nodes in few sampling times.
	This claim can be substantiated with a trivial induction.
	To this end, let \((\mathcal V,\mathcal E^t)\) denote the connectivity graph at time \(t\) among agents \(\mathcal V\), where \((j,i)\in\mathcal E^t\subseteq\mathcal V\times\mathcal V\) iff \(j\in\mathcal N_i^t\).
	The graph distance between agents \(i,j\in\mathcal V\) at time \(t\) is defined as
	\[
		\dist^{\mathcal E^t}(i,j)
	{}={}
		\min\set{\eta}[\exists(i,j_0),(j_0,j_1),\dots,(j_{\eta-1},j)\in\mathcal E^t],
	\]
	while the \emph{exact} hierarchy level of \(i\) at time \(t\) as
	\[
		D_i^{\mathcal E^t}
	{}={}
		\min\set{\overline\pi,\dist^{\mathcal E^t}(i,\ell)}[\ell\in\mathcal L],
	\]
	where \(\mathcal L\subseteq\mathcal V\) is the set of leader agents and \(\overline\pi\) is a chosen upper bound as in \eqref{eq:hierarchy}.

	\begin{lem}
		Suppose that the delay in exchanging information is bounded by \(\delta\geq0\) timesteps, and let \(D\in\set{0,\dots,\overline\pi}\).
		If \(\mathcal E^t\) has not changed in the last \((\delta+1)D\) steps before time \(\bar t\), then for every agent \(i\in\mathcal V\)
		\begin{equation}\label{eq:hierarchyDelay}
			\pi_i^{\bar t}~
			\begin{ifcases}
				~=D_i^{\mathcal E^{\bar t}} & D_i^{\mathcal E^{\bar t}}\leq D\\
				~>D \otherwise.
			\end{ifcases}
		\end{equation}
		\begin{proof}
			We proceed by induction on \(D\).
			If \(D=0\), then the claim holds trivially, since \(D^{\mathcal E^{\bar t}}=0\) holds iff \(i\) is a leader and is strictly greater otherwise.
			Suppose that the claim holds for \(D-1\in[0,\overline\pi-1]\), and suppose that \(\mathcal E^t\) remained constant in the last \((1+\delta)D\) steps before \(\bar t\).
			Let \(i\in\mathcal V\); without loss of generality, \(D_i^{\mathcal E^{\bar t}}\geq D\), for otherwise there is nothing to show.
			For \(j\in\mathcal N_i^{\bar t}\), let \(\delta_{j|i}\) be the delay with which \(i\) receives the information from \(j\) at time \(\bar t\), and observe that necessarily \(D_j^{\mathcal E^{\bar t}}\geq D-1\).
			Since \(\delta_{j|i}+1\leq(\delta+1)D\), it follows from the induction that
			\begin{equation}\label{eq:hierarchy-1}
				\pi_j^{\bar t-\delta_{j|i}-1}~
				\begin{ifcases}
					=D_j^{\mathcal E^{\bar t}}
				&
					D_j^{\mathcal E^{\bar t}}\leq D-1
					~~
					\text{\small(\(\Leftrightarrow D_j^{\mathcal E^{\bar t}}=D-1\))}
				\\
					\geq D
					\otherwise,
				\end{ifcases}
			\end{equation}
			where we used the assumption that \(\mathcal E^{\bar t-\delta_{j|i}-1}=\mathcal E^{\bar t}\).
			Then, by the update rule \eqref{eq:hierarchy},
			\begin{equation}\label{eq:hierarchyGeq}
				\pi_i^{\bar t}
			{}={}
				\min\set{
					\overline\pi,\,
					1+\min_{i\neq j\in\mathcal N_i^t}\pi_j^{\bar t-\delta_{j|i}-1}
				}
			{}\geq{}
				D.
			\end{equation}
			If \(D_i^{\mathcal E^{\bar t}}=D\), then there exists \(j\in\mathcal N_i^{\bar t}\) with \(D_j^{\mathcal E^{\bar t}}=D-1\).
			If instead \(D_i^{\mathcal E^{\bar t}}>D\), then \(D_j^{\mathcal E^{\bar t}}\geq D\) holds for all \(j\in\mathcal N_i^{\bar t}\).
			In both cases, \eqref{eq:hierarchyDelay} follows by combining \eqref{eq:hierarchy-1} and \eqref{eq:hierarchyGeq}.
		\end{proof}
	\end{lem}

	Note that including \(\overline\pi\) in the minimum also accounts for possible instances of isolatedness of agent \(i\) or detachedness of the cluster it belongs from a leader, while excluding \(\pi_i^{t-1}\) is necessary in order to guarantee an up-to-date (as opposed to a best-so-far) hierarchy estimate.
	The values of the own and neighbors' hierarchy levels are used to customize cohesion references, as described next.

\subsection{Subjective neighbor weighting and flocking references}\label{sec:weights}%
	In order to promote a flow towards the leader, each follower \(i\) can bias the reference center of mass dictating the cohesion rule based on the hierarchy level of all agents in its neighborhood, including itself.
	To this end, simplifying the notation as \(\pi_j\coloneqq\pi_i^t\) if \(j=i\), and \(\pi_j\coloneqq\pi_j^{t-1}\) otherwise, cohesion weights
	\begin{equation}\label{eq:w}
		w_{j|i}^{\tau|t}
	{}={}
		\frac{2^{-\pi_j}}{\sum_{\ell\in\mathcal N_i^{\tau|t}}2^{-\pi_\ell}}
	\end{equation}
	can be chosen to bias \(\bar p_i^{\tau|t}\) as in \eqref{eq:reference}.

	While the same can be done for the alignment weights \(\omega_{j|i}^{\tau|t}\), a more refined criterion can be implemented to discourage unwanted behaviors.
	For instance, in the Cucker-Smale model \cite{cucker2007emergent,shen2008cucker,zhang2016model}, \(\omega_{j|i}\) is inversely proportional to the distance between agents \(i\) and \(j\).
	While the rationale is perfectly reasonable and intuitive, when navigating an obstructed environment it is hardly helpful, as will be showcased in the simulations section.
	Differently from flocking, navigation presupposes a heading direction that is not the result of an agreement among agents, and a navigation-friendly alignment rule should be able to reflect this bias.
	Suppose that neighboring agent \(j\) is traveling \emph{ahead} of agent \(i\), and changes direction because of an obstacle in its way.
	In biasing its velocity accordingly, agent \(i\) is better prepared for avoiding the same obstacle in a later stage, even in case it is not currently detected.
	On the contrary, if agent \(j\) is \emph{behind}, then the obstacle is of no concern to agent \(i\), at least in the near future, and indulging into a possibly abrupt violation of the flow is not desirable in this case.
	Given that followers obey a flocking behavior and are unaware of the entire flock configuration, notions of \emph{ahead} and \emph{behind} are merely based on agents' instantaneous velocity.
	In this context, we may define subjective orientations by comparing the displacement between agents with the individual movement direction, and say that, at time \(t\) and relative to agent \(i\),
	\[
		j\in\mathcal N_i^{t|t}
		~
		\begin{ifcases}
			\text{is \emph{ahead}} & \innprod{v_i^{t|t}}{p_{j|i}^{t-1|t-1}-p_i^{t|t}}\geq0\\
			\text{is \emph{behind}} \otherwise.
		\end{ifcases}
	\]
	(The time difference between own and neighbor $j$'s positions $p_i^{t|t}$ and $p_{j|i}^{t-1|t-1}$ reflects the fact that data exchange among agents happens only once before the MPC problem is addressed; cf. workflow in \cref{alg:workflow}.)
	Once a behind weight \(\omega_{i,\rm b}\in[0,1]\) is fixed, we may define alignment weights as
	\begin{equation}\label{eq:alignment}
		\omega_{j|i}^t
	{}\coloneqq{}
		\begin{ifcases}
			1 & j \text{ is ahead wrt } i\\
			\omega_{i,\rm b} \otherwise
		\end{ifcases}
	\end{equation}
	to bias \(\bar v_i^{\tau|t}\) as in \eqref{eq:reference}.
	Our experiments suggest that values \(\omega_{i,\rm b}\leq 0.2\) tend to work well in all scenarios.

\subsection{Cohesion/alignment dynamic trade-off}\label{sec:tradeoff}%
	Once the reference positions and velocities are determined, it remains to decide if and how much tracking either quantity is to be favored.
	This is done by suitably choosing the weight matrix \(Q_i^t\) in the cost function \(J\) (recall that \(\bm y_i^t\) comprises position and velocity vectors).
	If agent \(i\) is far from the reference center of mass \(\bar p_i^{t|t}\), aligning with distant agents serves no purpose; in this case, \emph{cohesion} should be the priority reference.
	On the contrary, if agent \(i\) is very close to the center of mass, then cohesion is already well achieved and sticking to neighboring agents' \emph{alignment} becomes more reasonable.
	To dynamically adjust this trade-off, a Cucker-Smale-type weight rule \cite{cucker2007emergent} can be adopted, e.g.,
	\begin{equation}\label{eq:Q}
		\setlength\matrixcolsep{1pt}
		Q_i^t
	{}\coloneqq{}
		\text{\footnotesize\(\begin{pmatrix}1-q_i^t\\&1-q_i^t\\&&~q_i^t~\\&&&~q_i^t~\end{pmatrix}\)}
	~~\text{with}~~
		q_i^t
	{}\coloneqq{}
		\frac{q_{i,\rm st}}{1+c_i\|p_i^{t|t}-\bar p_i^{t|t}\|^2},
	\end{equation}
	for some \(c_i>0\) and \(q_{i,\rm st}\in(0,1)\); by setting \(c_i=0\) one recovers a conventional \emph{static} trade-off \(q_i^t\equiv q_{i,\rm st}\).
	Although with the same principle stage-dependent weights $Q_i^{\tau|t}$ could also be considered, experimental evidence has shown that this modification hardly makes any difference.
	This is somewhat expected, given that within short sampling times the changes in position and velocity are negligible.
	Note that, although based on the same distance inverse proportionality principle, this weighting rule is substantially different from the original Cucker-Smale criterion.
	In fact, while the latter is used to determine which velocity to align to, the proposed weight \eqref{eq:Q} quantifies the importance of said reference.

	\section{Simulations}\label{sec:simulations}%
		\begin{figure*}[t]
	\centering
	\renewcommand{\figname}{Roundabout}%
	\begin{addHFigTile}{}
		\caption{%
			\emph{Success of the proposed rules.}
			Despite the short detection ranges and the indifference of the leader, with the proposed rules all agents successfully manage to reach destination.%
		}%
	\end{addHFigTile}
	\hfill
	\begin{addHFigTile}{_NOdyn}
		\caption{%
			\emph{Failure of static trade-off weights \(Q_t\).}
			When far apart, reference velocities are less influential, while enforcing cohesion promotes proximity and consequent higher information.
		}%
	\end{addHFigTile}
	
	\vspace*{1pt}%
	\begin{addHFigTile}{_NOback_CSalign}
		\caption{%
			\emph{Failure of classic Cucker-Smale rule.}
			Velocity contributions of agents ``behind'' interfere with the ``ahead'' trajectory and should therefore weighed less.
			Orientation-blind weights easily lead to pathological ``drag-behind'' situations.
		}%
	\end{addHFigTile}
	\hfill
	\begin{addHFigTile}{_NOhier}
		\caption{%
			\emph{Failure of static hierarchies.}
			A sophisticated hiearchical chain as described in \eqref{eq:hierarchy} leads to a much faster spread of educated information, whereas prioritizing only leaders advantages only those directly connected.
			Far agents are then the culprit of a drag-behind tendency.%
		}%
		\label{fig:roundabout:NOhier}%
	\end{addHFigTile}
	\caption{%
		\emph{U-turn at a roundabout.}
		A leader travelling at a constant speed can be well tracked by all followers till destination.
		Neglecting some of the rules may result in disconnected agents left behind.
	}%
	\label{fig:roundabout}
\end{figure*}

\begin{figure}
	\centering
	\fbox{%
		\begin{minipage}{.97\linewidth}%
			\includegraphics[width=\linewidth]{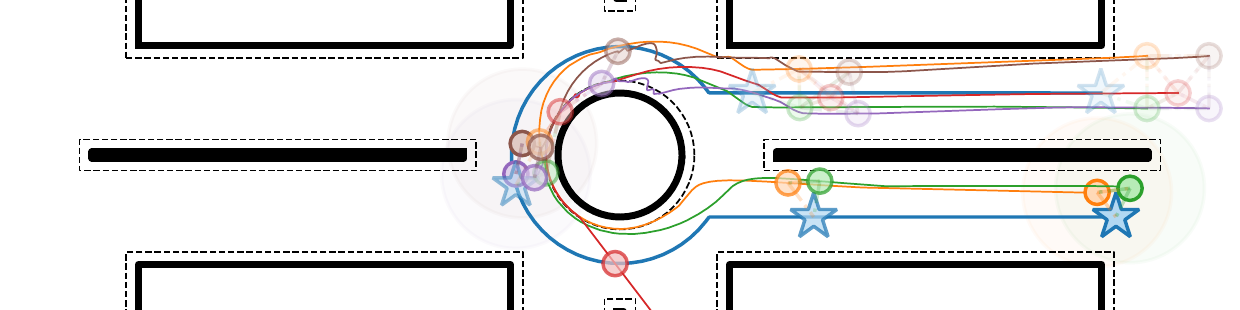}
		\end{minipage}%
	}%
	\caption{%
		\emph{Importance of multi-stage control}.
		The same roundabout scenario of \cref{fig:roundabout} fails without taking future predictions into account.
		In this example, the red follower loses all its neighbors and ends up pursuing a wrong direction.%
	}
	\label{fig:horizon}%
\end{figure}

To demonstrate the advantages of the purpose-oriented rules introduced in \Cref{sec:flocking}, we consider a system of planar point-mass boids with double integrator dynamics and one leader which navigate obstructed environments.
Except for the scenario depicted in \cref{fig:horizon}, which considers (the shortcomings of) a single-stage formulation, MPC horizons are set to \(T_i=8\) for all boids, while separation horizons as \(T_{i,\rm sep}=4\); avoidance of predicted collisions at stages higher than \(T_{i,\rm sep}\) is soft constrained with penalty \(\rho_{i,\rm sep}=100\).
The alignment weights of behind agents as in \eqref{eq:alignment} are set to \(\omega_{i,\rm b}=0.2\), while in the trade-off rule \eqref{eq:Q} we use a default static coefficient \(q_{i,\rm st}=0.5\) and a dynamic weight \(c_i=2\).
An upper cap \(\overline\pi=10\) is chosen for the hierarchy levels as in \eqref{eq:hierarchy}, \(\gamma_i=0.5\) is used as discount factor in the cost, and \(\Delta t=\frac{1}{40}~{\rm s}\) is set as sampling time.

The leader shares its predicted positions and velocities, but is otherwise indifferent to other agents.
We consider a proximity-based connectivity, i.e., \(j\in\mathcal N_i^t\) iff \(\|p_i^t-p_j^t\|\leq r\) for some detection radius \(r\), assumed the same for all agents for simplicity.
We consider bounded box constraints with radius 2 for (acceleration) inputs and velocity states, and obstacles of the form \(O=\set{x\in\R^2}[h^O(x)<0]\) for some smooth mapping \(\func{h^O}{\R^2}{\R^{n_O}}\).

\subsection{Implementation details}
	The experiments are run in Python, and the MPC problems are solved with the parametric optimizer OpEn \cite{sopasakis2020open} with direct interface.
	Following \cite{sathya2018embedded}, avoidance of obstacle \(O\) is modeled as the smooth and real-valued equation \(\left(\prod_{i=1}^{n_O}\min*\set{0,\,h^O_i}\right)^2=0\), and enforced with OpEn's penalty method, while velocity constraints and collision avoidance are enforced with the augmented Lagrangian method.
	Obstacles are slightly enlarged as a safety precaution; the enlargements are represented by dashed lines in the figures.
	Parameters passed to the solver are initial state (position and velocity) \(s_i^{t|t}\in\R^4\), cohesion and alignment references \(\bar{\bm p}_i^t,\bar{\bm v}_i^t\in\R^{2T}\), predicted positions of the \(J=5\) closest neighbors (excluding \(i\) itself) \(\bm p_{j|i}^t\in\R^{2T}\), and cohesion/alignment trade-off weight \(q_i^t\), which we saturate in \([0.2,0.8]\) as numerical evidence suggests that too extreme trade-offs may lead to erratic behaviors.
	The limit on the number of neighbors to separate from is imposed due to modeling restrictions, in account of the non variable size of the parameter vector.
	For the same reason, when \(\mathcal N_i^{\tau|t}\setminus\set i\) has less than \(J\) agents, parameters are filled with \(10^4\) distance units for dummy positions.

\begin{figure*}[t]
	\centering
	\renewcommand{\figname}{S_shape}%
	\begin{addVFigTile}{}
		\caption{%
			\emph{Success of the proposed rules.}
			All agents seamlessly manage to navigate the impervious environment until destination.%
		}%
	\end{addVFigTile}
	\hfill
	\begin{addVFigTile}{_NOdyn}
		\caption{%
			\emph{Failure of static trade-off weights \(Q_t\).}
			At the second obstacle, the flock looses momentum by not prioritizing cohesion.%
		}%
	\end{addVFigTile}
	\hfill
	\begin{addVFigTile}{_NOback_CSalign}
		\caption{%
			\emph{Failure of classic Cucker-Smale rule.}
			``\(\alpha\)-lattice'' formations \cite{olfatisaber2006flocking} resulting from traditional flocking are not fit for obstructed environments.%
		}%
	\end{addVFigTile}
	\hfill
	\begin{addVFigTile}{_NOhier}
		\caption{%
			\emph{Failure of static hierarchies.}
			Similarly to the situation in \figref{fig:roundabout:NOhier}, trailing agents drag the flock behind, unless a smart hiearchy is established.%
		}%
	\end{addVFigTile}
	\caption{%
		\emph{Tracking an inaccurate virtual reference through an obstructed narrow passage.}
		The leader is here a virtual reference which is not accurately detected or computed by the followers.
		As a result, it violates the first obstacle on the way.
	}
	\label{fig:S_shape}
\end{figure*}

\subsection{Problems description}
	We consider two scenarios that are hostile to conventional flocking control problems; one is taking a U-turn at a roundabout, and the other is going through an obstructed narrow passage.
	In these situations, the traditional flocking rules may not be meaningful or even disadvantageous, as they discourage changes in the formation which are instead essential in navigation:
	a U-turn imposes a quick change of direction subject to road constraints, i.e., lanes orientation and boundaries,
	while obstructed environments make it easier for boids to get stuck by an obstacle and then either dragging back the rest of the flock or ending up disconnected from it.
	Instead, in the sequel we demonstrate that the proposed rules well cope with these scenarios, even when tracking an imprecise leading reference.
	Indeed, when any of the rules is inactive, the flock fails to reach the destination, while the complete method successfully achieves the task in both cases.
	
	The leader agent is denoted with a blue star, while the followers with colored circles.
	The leader proceeds along a predefined path at a velocity compatible with the followers' limits, but is otherwise indifferent to their states.
	Snapshots of the configurations at different time instants are shown, where faded colors indicate the ealier time.
	Agents are assumed dimensionless point particles positioned in the center of the circles, the radiuses of which indicate the separation distance \(d_{i,\rm sep}\).
	Dashed lines between agents indicate the (time-varying, distance-based) connectivity graph; the communication ranges, depicted as transparent circles in the last snapshot, are set very small to make the problem more challenging.
	Obstacles are represented by a combination of rectangles and circles, shown in black, with respective enlargements as dotted lines.
	Results of the simulations are described in the captions of the respective \cref{fig:horizon,,fig:roundabout,,fig:S_shape}.

	Only few agents are considered for the sake of graphical clarity.
	We however emphasize that the advantage of the proposed rules is more evident when considering large fleets, which are more prone to exhibiting disconnections in the communication graph.

	\section{Conclusions}\label{sec:conclusions}%
		A distributed control strategy using MPC has been proposed for a flock navigation problem in unknown obstructed environments.
Conventional flocking rules have been suitably modified, thus introducing a dynamic hierarchy strategy, subjective neighbor weighting, and cohesion/alignment dynamic trade-off that are all computed using local information only.
Numerical examples demonstrate that these rules effectively guide a flock in challenging environments, and without imposing severe restrictions on leaders' trajectories.
Possible future directions include accounting for noisy/faulty communication and investigating optimal cooperative leaders' control actions.
It would also be interesting to extend the analysis of (static) \emph{layered-path graphs} of \cite{yoshise2022algebraic} to our dynamic hierarchy setting, so as to develop strategies to minimize connectivity loss in case of agents leaving the flock.%

	\bibliographystyle{plain}%
	\bibliography{TeX/Bibliography.bib}

\end{document}